\documentclass[11pt, reqno]{amsart}

\usepackage{amssymb,latexsym,amsmath,amsfonts,amsthm}
\usepackage{mathrsfs}
\usepackage{enumitem}
\usepackage[usenames]{color}
\usepackage{hyperref}
\usepackage{comment}

\allowdisplaybreaks

\numberwithin{equation}{section}

\definecolor{DPurple}{rgb}{0.46,0.2,0.69}

\theoremstyle{definition}
\newtheorem{definition}{Definition}[section]

\theoremstyle{remark}
\newtheorem{remark}[definition]{Remark}

\theoremstyle{plain}
\newtheorem{theorem}[definition]{Theorem}
\newtheorem{result}[definition]{Result}
\newtheorem{lemma}[definition]{Lemma}
\newtheorem{proposition}[definition]{Proposition}

\newcommand{\zt}{\zeta}

\newcommand{\poin}{\boldsymbol{{\sf p}}}

\newcommand{\D}{\mathbb{D}}

\newcommand{\smoo}{\mathcal{C}}

\newcommand{\lrarw}{\longrightarrow}

\newcommand{\bdy}{\partial}

\newcommand{\Cn}{\mathbb{C}^n}
\newcommand{\CC}{\mathbb{C}^2}
\newcommand{\C}{\mathbb{C}} 
\newcommand{\R}{\mathbb{R}}

\newcommand{\wt}{\widetilde}

\begin{document}

\title[Incomplete taut visibility domains]{Taut visibility domains are not necessarily \\ Kobayashi complete}

\author{Rumpa Masanta}
\address{Department of Mathematics, Indian Institute of Science, Bangalore 560012, India}
\email{rumpamasanta@iisc.ac.in}

\begin{abstract}
We answer a question asked recently by Banik in the negative by showing that for each $n\geq 2$, there exists
a taut visibility domain in $\Cn$ that is not Kobayashi complete. The domains that we produce are bounded and
have boundaries that are very regular away from a single point.
\end{abstract}


\maketitle

\vspace{-4mm}
\section{Introduction}\label{S:intro}
The word ``visibility'' in our title refers to a notion introduced by Bharali--Zimmer
\cite{bharalizimmer:gdwnv17, bharalizimmer:gdwnv23}\,---\,also see \cite{bharalimaitra:awnovfoewt21} by 
Bharali--Maitra\,---\,and which has generated a lot of interest in exploring the complex geometry of domains
using the Kobayashi distance. The visibility property can be seen as a weak notion of negative curvature as
it resembles visibility in the sense of Eberlein--O'Neill \cite{eberleinneill:vm73}. Intuitively, the
visibility property requires that all geodesics with end-points close to two distinct points in the boundary
must bend uniformly into the domain. A precise definition of the above property will be given below. This
property has many applications where one needs to control the behaviour of certain classes of holomorphic maps into
domains with the visibility property. These applications are now too numerous to list in this short note.
The property itself, i.e., its general consequences, has been investigated in a number of recent papers: see, 
for instance, 
\cite{bharalizimmer:gdwnv17, bharalimaitra:awnovfoewt21, braccinikolovthomas:vkgcdrp22, bharalizimmer:gdwnv23,
chandelmitrasarkar:nvwrkdca21, masanta:vdekdcm24}.  
\smallskip

For a domain $\Omega\varsubsetneq\Cn$, it is in general very difficult, when $n\geq 2$, to determine whether
the domain $\Omega$ equipped with the Kobayashi distance $K_\Omega$ is Cauchy-complete\,---\,and, so,
whether $(\Omega,K_\Omega)$ is a geodesic space. Thus, the formal
definition of visibility property needs some care. We shall give this definition after stating our result.
Dealing with the issue of whether $(\Omega,K_\Omega)$ is Cauchy-complete often leads to proofs that can be
quite technical. This motivated Banik, in a recent article \cite{banik:vdtanp24} to ask 
(more on this below): 
\emph{Let $\Omega\varsubsetneq\Cn$, $n\geq 2$, be a bounded visibility domain. Assume that $\Omega$ is taut.
Then, is $\Omega$ Kobayashi complete?} We answer this question in the negative, as follows:

\begin{theorem}\label{T:main-result}
For each $n\geq 2$, there exists a bounded domain $\Omega\varsubsetneq\Cn$ with the properties
\begin{itemize}
  \item $0\in \bdy\Omega$ and $\bdy\Omega$ is $\smoo^\infty$-smooth away from $0$,
  \item $\bdy\Omega$ is strongly Levi-pseudoconvex at each point $p\in \bdy\Omega\setminus\{0\}$,
\end{itemize}
such that $\Omega$ is taut and is a visibility domain with respect to the Kobayashi distance, but
$(\Omega,K_\Omega)$ is not Cauchy-complete.  
\end{theorem}

To understand what motivates the question above, we need a couple of definitions. Since our aim is to
emphasise the \textbf{ideas} involved, our second definition will be restricted to bounded domains in $\Cn$.
\smallskip

Given a domain $\Omega\varsubsetneq\Cn$, we say that $\Omega$ is \emph{Kobayashi hyperbolic} if the Kobayashi
pseudodistance $K_\Omega$ is a distance. Let $k_\Omega:\Omega\times\Cn\lrarw[0,\infty)$ denote the
Kobayashi pseudometric on $\Omega$. With these, we can now give the following definitions. 

\begin{definition}\label{D:almost-geodesic}
Let $\Omega\varsubsetneq\Cn$ be a Kobayashi hyperbolic
domain. Let
$I\subseteq \R$ be an interval. For
$\lambda\geq 1$ and $\kappa\geq 0$, a curve $\sigma: I\lrarw \Omega$ is said to be a 
\emph{$(\lambda, \kappa)$-almost-geodesic} if 
\begin{itemize}[leftmargin=25pt]
 \item[$(a)$] for all $s,t\in I$
 \[
   \frac{1}{\lambda}|s-t|-\kappa\leq K_\Omega (\sigma(s),\sigma(t))\leq \lambda|s-t|+\kappa,
 \]
 \item[$(b)$] $\sigma$ is absolutely continuous (whence $\sigma'(t)$ exists for almost every $t\in I$),
 and for almost every $t\in I$,
 $k_\Omega(\sigma(t);\sigma'(t))\leq\lambda$.
\end{itemize}
\end{definition}

\begin{definition}\label{D:visibility-domain}
Let $\Omega\varsubsetneq\Cn$ be a bounded domain. We say that 
\emph{$\Omega$ is a visibility domain with respect to the Kobayashi distance} (or simply a 
\emph{visibility domain} for brevity) if, for any
$\lambda\geq 1$, $\kappa\geq 0$, each pair of distinct points $p,q\in\bdy\Omega$, and each pair of
$\overline\Omega$-open neighbourhoods
$U_p$ of $p$ and $U_q$ of $q$ such that $\overline{U_p}\cap\overline{U_q}=\emptyset$, there exists a compact set
$K\subset\Omega$ such that
any $(\lambda,\kappa)$-almost-geodesic $\sigma:[0,T]\lrarw\Omega$ with $\sigma(0)\in U_p$
and $\sigma(T)\in U_q$, we have $\sigma([0,T])\cap K\neq\emptyset$.
\end{definition}

We already mentioned above the difficulty of knowing whether $(\Omega,K_\Omega)$ is Cauchy-complete\,---\,and,
consequently, whether $(\Omega,K_\Omega)$ is a geodesic space. Therefore, geodesics joining 
any two distinct points in $\Omega$ 
may not exist. Thus, $(\lambda,\kappa)$-almost-geodesics serve as substitutes for the role of geodesics
because of this result by Bharali--Zimmer \cite[Proposition~5.3]{bharalizimmer:gdwnv23}: \emph{if 
$\Omega\varsubsetneq\Cn$ is a Kobayashi hyperbolic domain, then for
any $\kappa> 0$, every two points in $\Omega$ are joined by a (1, $\kappa$)-almost-geodesic.}
\smallskip

Now that we have seen that a $(\lambda,\kappa)$-almost-geodesic is a considerably more technical object than
a geodesic, but that working with the former is unavoidable, one can, perhaps, appreciate why proofs 
involving $(\lambda,\kappa)$-almost-geodesics can be analytically complicated. Thus, the question, 
\emph{``Are all visibility domains Kobayashi complete?''} has come up in the past. Banik's results answer
this question in the negative, which calls for a refinement of the above question. We refer the reader to
\cite[Section~5]{banik:vdtanp24} for Banik's motivations for raising two separate questions. 
Theorem~\ref{T:main-result} answers \cite[Question~5.2]{banik:vdtanp24}. These questions appear to have
raised a lot of interest: very recently, a partial answer to \cite[Question~5.1]{banik:vdtanp24} was given
by Nikolov \emph{et al.} \cite{nikolovoktenthomas:vsdp24}.
\medskip

\section{Preliminaries}
This section is devoted to a few facts that will be needed in proving our main theorem. 
But we first clarify
some notation: for $z\in\Cn$, $\|z\|$ denotes the Euclidean norm on $\Cn$, and $\mathbb{B}^n(z,r)$ denotes
the open Euclidean
ball centered at $z$ with radius $r>0$.

\begin{lemma}\label{L:kobayashi-identity}
Let $G$ be a domain in $\C^m$ and $\Omega$ be a domain in $\Cn$, where $n>m$, such that 
$G\times\{0\}\subset \Omega\subset G\times\C^{n-m}$. Then, for any $z,w\in G$, 
$K_\Omega((z,0),(w,0))=K_G(z,w)$.
\end{lemma}

The proof for Lemma~\ref{L:kobayashi-identity} follows from the contractive property for the Kobayashi 
pseudodistance of the obvious inclusion 
$G\hookrightarrow\Omega$ and the product property. 
\smallskip

For the next result, we need a very specific construction. Write:
\[
  a_\nu:=1/{4^{\nu+1}}, \quad
  b_\nu:=1/2^{\nu+1},
\]
and define
\[
  X_{\nu} := \{(z_1,z_2)\in\CC : z_2 = (a_{\nu+1}+a_\nu)z_1 -a_\nu a_{\nu+1}, \ |z_1|\leq b_\nu\},
\]
$\nu = 1,2,3,\dots$, and let $X := \bigcup_{\nu\geq 1}X_\nu$. Consider the holomorphic maps 
$\phi_\nu: \D\lrarw \C^2$, 
defined as $\phi_\nu(\zt) := (b_\nu\zt,\,b_\nu(a_{\nu+1}+a_\nu)\zt - a_\nu a_{\nu+1})$ for each $\zt\in \D$.
If $x_\nu := (a_\nu, a_{\nu}^2)$, then the objects defined have the following properties:
\begin{itemize}
  \item[$(a)$] $(x_\nu)_{\nu\geq 1}\subset X$ and $x_\nu\to 0$ as $\nu\to \infty$, \smallskip
  
  \item[$(b)$] $\phi_\nu(\D)\subset X$ for $\nu = 1,2,3,\dots$, \smallskip
  
  \item[$(c)$] $a_\nu/b_\nu, a_{\nu+1}/b_\nu \in \D$, $\phi_\nu(a_\nu/b_\nu) = x_\nu$, and
  $\phi_\nu(a_{\nu+1}/b_\nu) = x_{\nu+1}$ for $\nu = 1,2,3,\dots$, and \smallskip
  \item[$(d)$] $\sum_{\nu=1}^\infty \poin(a_\nu/b_\nu,a_{\nu+1}/b_{\nu}) < \infty$,
\end{itemize}
where $\poin$ denotes the Poincar{\'e} distance.

\begin{proposition}\label{P:psh-constr-c2}
There exists a continuous, proper plurisubharmonic function $u: \mathbb{B}^2(0,3)\lrarw \R$ with the properties
\begin{itemize}
  \item[$(a)$] $u$ is $\smoo^\infty$-smooth and strictly plurisubharmonic on $\mathbb{B}^2(0,3)\setminus\{0\}$,
  \item[$(b)$] $\nabla u(z)\neq 0$ for $z\in\mathbb{B}^2(0,3)\setminus\{0\}$, $u(0)=1$, and
  \item[$(c)$] $u|_{X} < 1$ (where $X$ is the set described above),
\end{itemize}
such that if $G\varsubsetneq\C^2$ denotes the connected component of $\{z\in\mathbb{B}^2(0,3):u(z)<1\}$
containing the set $X$,
then $G\Subset \mathbb{B}^2(0,3)$, $G$ is pseudoconvex, and such that the sequence 
$(x_\nu)_{\nu\geq 1}$ (described above)
is $K_{G}$-Cauchy.
\end{proposition}

The argument for Proposition~\ref{P:psh-constr-c2} is as in the proof of
\cite[Theorem~7.5.9]{jarnickipflug:idmca93} (also see \cite{jarnickipflug:cKcbd91}), which the above paraphrases. The latter result,
as is
mentioned in \cite{jarnickipflug:idmca93}, is an unpublished result of Nessim Sibony.%
\smallskip

We shall need the notion of \emph{local Goldilocks points}, introduced in \cite{bharalizimmer:gdwnv23}, in
order to prove Theorem~\ref{T:main-result}. However, our proof will reference a very special type of local 
Goldilocks point. Thus, in the interests of brevity, we shall not define the above notion, but refer the 
reader to \cite[Section~1.1]{bharalizimmer:gdwnv23}. Instead, we will state a result of Bharali--Zimmer 
\cite{bharalizimmer:gdwnv23} (where they denote the set of all local Goldilocks points in $\bdy\Omega$
by $\bdy_{{\rm lg}}\Omega$), which is relevant for our proof.

\begin{result}[paraphrasing {\cite[Theorem~1.4]{bharalizimmer:gdwnv23}}]\label{R:local-goldilocks}
Let $\Omega\varsubsetneq\Cn$ be a bounded domain. Suppose the set $\bdy\Omega\setminus\bdy_{{\rm lg}}\Omega$ is
totally disconnected. Then, $\Omega$ is a visibility domain with respect to the Kobayashi distance.
\end{result}
\smallskip

\section{The proof of Theorem~\ref{T:main-result}}
We now present the proof of our result.

\begin{proof}[The proof of Theorem~\ref{T:main-result}]
Let $u$, $G$, $X$, and $(x_\nu)_{\nu\geq 1}$ be as in Proposition~\ref{P:psh-constr-c2}. First, we will
prove the result for $n\geq 3$. Define a
function $h:\mathbb{B}^2(0,3)\times\C^{n-2}\lrarw\R$ by 
\[
  h(z_1,z_2,...,z_n):=u(z_1,z_2)+\sum_{j=3}^{n}|z_j|^2.
\]
Let us define the set $\mathcal{U}:=(\mathbb{B}^2(0,3)\setminus\{0\})\times\C^{n-2}$. Clearly, $h$ is a
continuous, plurisubharmonic function (since $u$ is so). 
Since $u$ satisfies the properties $(a)$ and $(b)$ in
Proposition~\ref{P:psh-constr-c2}, it follows that $h$ is 
$\smoo^\infty$-smooth and strictly plurisubharmonic on $\mathcal{U}$, 
$\nabla h(z)\neq 0$ for each $z\in\mathcal{U}$, and $h(0)=1$.  Let 
$S:=X\times\{0\}\subset\Cn$. Since
$u{\mid}_X<1$, $h{\mid}_{S}<1$. Let $\Omega$ be the connected component of the open set 
$\{z\in\mathbb{B}^2(0,3)\times\C^{n-2}:h(z)<1\}$ that contains the connected set $S$. Clearly, 
$G\times\{0\}\subset\Omega$. If $p\in\bdy\Omega$, then $p\in\mathcal{U}$ unless $p=(0,...,0)$. 
Thus, as $u$ is proper, $\nabla h(p)$ is defined for $p\in\bdy\Omega\setminus\{0\}$ and, from the above
properties of $h$,
it follows that $\Omega$ is bounded and $\bdy\Omega$ is $\smoo^\infty$-smooth away from $0$.
Similarly, since $h$ is strictly plurisubharmonic at each point $z\in\mathcal{U}$, we have
$\bdy\Omega$ is strongly Levi-pseudoconvex at each point $p\in\bdy\Omega\setminus\{0\}$. By the construction
of $h$, we have 
$G\times\{0\}\subset \Omega\subset G\times\C^{n-m}$.
We shall now complete our argument in 3 steps.
\medskip

\noindent{{\textbf{Step~1.}}}
\emph{$\Omega$ is taut.}
\smallskip

\noindent{Since $\Omega$ admits a continuous plurisubharmonic exhaustion function $h$, which is bounded
above on $\Omega$, it follows that $\Omega$
is a hyperconvex domain. Hence, it is taut.}
\medskip

\noindent{{\textbf{Step~2.}}}
\emph{$\Omega$ is a visibility domain.}
\smallskip

\noindent{Since $\Omega$ is bounded, $\Omega$ is a Kobayashi hyperbolic domain. From the discussion above, we
have $\bdy\Omega$ is
strongly Levi-pseudoconvex at each point $p\in\bdy\Omega\setminus\{0\}$. It has been shown
in the proof of \cite[Theorem~1.2]{banik:vdtanp24} that every strongly Levi-pseudoconvex point is a local
Goldilocks point. Thus,
every $p\in\bdy\Omega\setminus\{0\}$ is a local Goldilocks point. Therefore, by Result~\ref{R:local-goldilocks},
$\Omega$ is a visibility domain.}%
\medskip

\noindent{{\textbf{Step~3.}}}
\emph{$(\Omega,K_\Omega)$ is not Cauchy-complete.}
\smallskip

\noindent{Write $\wt x_\nu:=(x_\nu,0,...,0)\in S$ for each $\nu\geq 1$. Since $S\subset\Omega$, 
$(\wt x_\nu)_{\nu\geq 1}\subset\Omega$. By Lemma~\ref{L:kobayashi-identity}, any sequence
$(z_\nu)_{\nu\geq 1}\subset G$ is $K_G$-Cauchy if and only if $((z_\nu,0,...,0))_{\nu\geq 1}\subset\Omega$ is
$K_\Omega$-Cauchy. Thus, by Proposition~\ref{P:psh-constr-c2}, $(\wt x_\nu)_{\nu\geq 1}$ is 
$K_\Omega$-Cauchy. However, it is not convergent in $\Omega$ as $\wt x_\nu\to 0\in\bdy\Omega$. Hence,
$(\Omega,K_\Omega)$ is not Cauchy-complete.}
\smallskip

Therefore, $\Omega$ is the desired domain, and the result follows for $n\geq 3$.
\smallskip

For $n=2$, by Proposition~\ref{P:psh-constr-c2}, we can follow the same argument in Steps~1--3 by
replacing $h$ by $u$, $\Omega$ by $G$, and $S$ by $X$. 
Hence the result.
\end{proof}

\begin{remark}
The reader might notice that we do not really need Lemma~\ref{L:kobayashi-identity} to carry out the
argument in Step~3 above. However, it is amusing, and allows one to justify Step~3 without using any
inequalities.  
\end{remark}
\smallskip

\section*{Acknowledgements}
I would like to thank my thesis advisor, Prof. Gautam Bharali, 
for his helpful observations on shortening the exposition above. This work is
supported in part by a scholarship from the Prime
Minister's Research Fellowship (PMRF) programme (fellowship no.~0201077) 
and by a DST-FIST grant (grant no.~DST FIST-2021 [TPN-700661]).

\end{document}